\documentclass[a4paper, 11pt]{article}

\usepackage{amsmath}
\usepackage{amssymb}
\usepackage{graphics}
\usepackage[dvips]{graphicx}
\usepackage{epsfig}
\usepackage[english]{babel}
\usepackage{lineno}
\usepackage{color}
\usepackage{caption}
\usepackage{subcaption}
\usepackage{epstopdf}

\begin{document}

\title{Super compact pairwise model for SIS epidemic on heterogeneous networks}

\author{P\'eter L. Simon $^{1,\ast}$, Istvan Z. Kiss $^{2}$}

\maketitle

\begin{center}
$^1$ Institute of Mathematics, E\"otv\"os Lor{\'a}nd University Budapest, and \\ Numerical Analysis and Large Networks Research Group, Hungarian Academy of Sciences, Hungary\\
$^2$ School of Mathematical and Physical Sciences, Department of Mathematics, University of Sussex, Falmer, Brighton BN1 9QH, UK
\end{center}

\vspace{1cm}

\begin{abstract}
In this paper we provide the derivation of a super compact pairwise model with only $4$ equations in the context of
describing susceptible-infected-susceptible ($SIS$) epidemic dynamics on heterogenous networks.
The super compact model is based on a new closure relation that involves not only the average degree but also the second and third moments
of the degree distribution. Its derivation uses an a priori approximation of the degree distribution of susceptible nodes in
terms of the degree distribution of the network. The new closure gives excellent agreement with heterogeneous
pairwise models that contain significantly more differential equations.\\
\end{abstract}

\noindent {\bf Keywords:} SIS epidemic; pairwise model, triple closure \\

\vspace{1cm}


\vspace{1cm}
\begin{flushleft}
$\ast$ corresponding author\\
email: simonp@cs.elte.hu\\
\end{flushleft}

\newpage

\section{Introduction}

While networks have provided a new modelling paradigm for population dynamics \cite{keeling2005networks,danon2011networks,pastor2014epidemic}, these are still used in
conjunction with mean-field models of various types. The most frequently used and well-known mean-field models for network epidemics are the degree-based mean-field (DBMF)
model, also known as heterogeneous mean-field \cite{SatorrasVespignani,pastor2014epidemic} and pairwise model \cite{rand2009,Keeling1999,House2011}. Both continue to provide a
productive framework for approximating expected values of random variables emerging from explicit network-based stochastic simulations in different contexts and networks with
different properites. The major advantage of such mean-field models stems from the fact that often these allow us to analytically determine quantities such
as the basic reproduction number, final epidemic size or endemic equilibrium \cite{SatorrasVespignani,Keeling1999}. Such analytic expressions then lead to a significantly better understanding of the interplay between network and disease characteristics.

Pairwise models have originally been introduced in the context of mathematical ecology \cite{sato1994pathogen} followed by natural extensions to
epidemiology \cite{Keeling1999}. The original simple model for undirected and unweighted networks has been subsequently extended to networks with heterogenous degree
\cite{eames2002modeling}, directed networks \cite{sharkey2006pair}, weighted  networks \cite{rattana2013class}, networks displaying motifs \cite{house2009motif} and even
combined with the edge-based compartmental modelling framework for an even more compact treatment \cite{House2011}.

\begin{figure}[h!]
  \centering
	\includegraphics[width=0.6\textwidth, height=0.3\textheight]{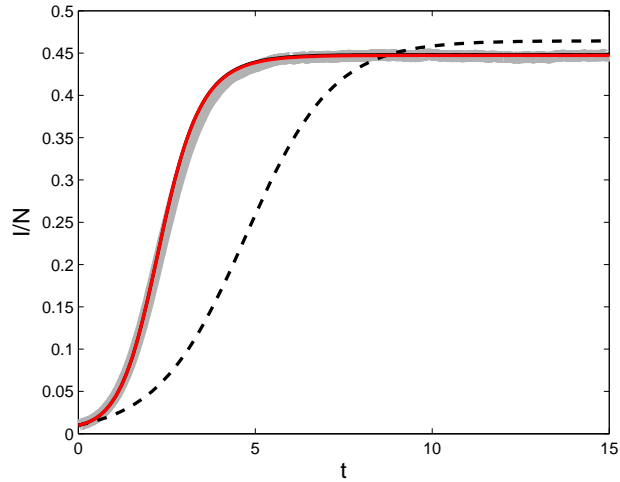}
	\caption{SIS epidemic propagation on a bimodal configuration random graph: simulation (gray thick curve), pair-wise (black dashed), compact pair-wise (black continuous), heterogeneous pair-wise (red continuous). The two latter curves are nearly indistinguishable. The parameter values are $N=1000$, $N_1=N_2=N/2$, $k_1=5$, $k_2=35$, $\gamma=1$ and $\tau = 3 \gamma \langle k \rangle/\langle k^2 \rangle$.} \label{fig:fig1}
\end{figure}

The closure in the most basic or fundamental pairwise model is based on the assumption on homogeneity of the degree distribution, i.e., all nodes have approximately the same
degree $n$. Hence the traditional pairwise model cannot be applied for graphs with heterogeneous degree distribution, such as bimodal graphs or networks with power law degree
distribution. This is shown in Figure \ref{fig:fig1}. For heterogeneous networks, a corresponding pairwise model was introduced in \cite{eames2002modeling}. This gives excellent
agreement with simulations for all configuration-like random networks \cite{molloy1995critical}, see Figure \ref{fig:fig1}. The heterogeneous pairwise model consist of order $N^2$
differential equations, where $N$ denotes the number of nodes in the network. An approximation of pairs leads to a simpler system, called compact pair-wise model that consist of only order $N$ equations \cite{House2011} and still gives very good agreement with simulations, see Figure \ref{fig:fig1}.

The aim of this paper is to introduce an even simpler model with only four equations that performs well for large heterogeneous networks.
The system is derived from the compact pairwise model by introducing a further approximation, and using a closure
relation that contains not only the average of the network's degree distribution but also its second and third moments.

\section{Derivation of the super compact paiwise (PW) model}

\subsection{Pairwise model for homogenous networks}
We start from the exact PW model. For the $SIS$ epidemic on an arbitrary undirected network the expected values of $[S]$, $[I]$, $[SI]$, $[II]$ and $[SS]$ satisfy the following system of differential equations
\begin{eqnarray}
\dot {[S]} &=& \gamma [I] - \tau [SI], \label{SIpairde1}\\
\dot {[I]} &=& \tau [SI]- \gamma [I], \label{SIpairde2}\\
\dot {[SI]} &=& \gamma ([II]-[SI]) + \tau ([SSI]-[ISI]-[SI]), \label{SIpairde3}\\
\dot {[SS]} &=& 2\gamma [SI] -2 \tau [SSI] \label{SIpairde4},\\
\dot {[II]} &=& -2\gamma [II] + 2\tau ([ISI]+[SI]), \label{SIpairde5}
\end{eqnarray}
where $[X]$, $[XY]$ and $[XYZ]$ denote the expected number of nodes in state $X$, edges in state $X-Y$ and triples in state $X-Y-Z$, respectively, with counting
according to all possible edge directions. This effectively means that for undirected networks $[XY]=[YX]$, $[XX]$ is double the number of unique edges in state $X-X$, and
similarly $X-Y-X$ accounts twice for a unique $X-Y-X$ triple, where $X,Y\in \{S,I\}$. This system is derived directly from master equations in \cite{taylor2012markovian} and hence exact. We note that some of
the equations can be omitted by exploiting conservation identities, such as $[S]+[I]=N$.

It is well known that in order to transform Eqs.~\eqref{SIpairde1}-\eqref{SIpairde5} into a self consistent solvable system closures need to be applied in order to break dependency on
higher order moments.  Particularly useful are closures at the level of triples. As it is well-known, the simplest closure is
\begin{equation}
[ASI]\approx \frac{n-1}{n} \frac{[AS][SI]}{[S]} ,
\label{triple_closure1}
\end{equation}
where $n=\langle k \rangle$ is the average degree of the network, and $A$ stands for $S$ or $I$. This closure leads to the traditional pairwise system
\begin{eqnarray}
\dot {[S]_p} &=& \gamma [I]_p - \tau [SI]_p, \label{SIpairclosde1}\\
\dot {[I]_p} &=& \tau [SI]_p- \gamma [I]_p, \label{SIpairclosde2}\\
\dot {[SI]_p} &=& \gamma ([II]_p-[SI]_p) + \tau \frac{n-1}{n}  \frac{[SI]_p ([SS]_p-[SI]_p)}{[S]_p} -\tau [SI]_p, \label{SIpairclosde3}\\
\dot {[SS]_p} &=& 2\gamma [SI]_p -2 \tau \frac{n-1}{n} \frac{[SI]_p [SS]_p}{[S]_p}, \label{SIpairclosde4} \\
\dot {[II]_p} &=& -2\gamma [II]_p + 2\tau \frac{n-1}{n} \frac{[SI]_p^2}{[S]_p}+ 2\tau [SI]_p . \label{SIpairclosde5}
\end{eqnarray}
Here the subscript $p$ is used to emphasize that the solution of this system is different from the exact values of the expected variables. As Figure 1 shows, this system cannot capture
network heterogeneities, hence closure \eqref{triple_closure1} needs improvement.

\subsection{Pairwise models for heterogenous networks: the heterogeneous, pre-compact and compact pairwise models}

The problem with closure \eqref{triple_closure1} is that it assumes that each node has degree $n$, which is obviously a crude approximation for heterogeneous networks. This has led to several heterogeneous mean-field models, where the state space is much extended to account for the expected number of nodes in different states and with a given degree, e.g., $[S_k](t)$ and $[I_k](t)$ for the expected number of susceptible and infected nodes of degree $k$, respectively. These new variables will induce or require further variables at pair level, such as $[S_kI_l](t)$ which denotes the expected value of the number of edges connecting susceptible nodes of degree $k$ to infected nodes of degree $l$. In this spirit, the following heterogeneous models were developed in historical order:
\begin{itemize}
\item heterogeneous pairwise model \cite{eames2002modeling},
\item pre-compact pairwise model \cite{eames2002modeling} and
\item compact pairwise model \cite{House2011}.
\end{itemize}

Instead of presenting the systems of differential equations of these models and working from the most explicit or complex to the more compact one, we start from the simplest model and show in an intuitive way how the more sophisticated models arise. Since closure \eqref{triple_closure1} uses the degree of the middle node, it is useful to express the triple as
$$
[ASI] = \sum_{k=1}^K [AS_kI],
$$
where the different degrees occurring in the graph are $k=1,2, \ldots , K$. The closure for the the triples in the right hand side can be written as
\begin{equation}
[AS_kI]\approx \frac{k-1}{k} \frac{[AS_k][S_kI]}{[S_k]} .
\label{triple_closurek}
\end{equation}
In order to use this closure in the exact system \eqref{SIpairde1}-\eqref{SIpairde5} one needs differential equations for $[S_k]$, for $[S_kI]$ and for $[S_kS]$. The exact
differential equations for $[S_k]$ are
\begin{eqnarray}
\dot {[S_k]} &=& \gamma [I_k] - \tau [S_kI], \quad k=1,2, \ldots , K, \label{SIhetde1}
\end{eqnarray}
where the substitution $[I_k] = N_k-[S_k]$ can be used. The simplest heterogeneous model \cite{House2011} uses only $[S_k]$ as new variables and introduces an algebraic expression that approximates $[S_kI]$ and $[S_kS]$ in terms of $[S_k]$, $[SI]$ and $[SS]$ as follows:
\begin{equation}
[S_kI] \approx [SI] \frac{k[S_k]}{\sum_{l=1}^K l[S_l]},
\label{cPWappr}
\end{equation}
which can be interpreted as showing that the ratio of the number of edges connecting degree $k$ susceptible nodes to infected nodes and the number of $SI$ edges is almost the same as the ratio of the number of stubs starting from degree $k$ susceptible nodes and the total number of stubs starting from susceptible nodes. Using this approximation, closure \eqref{triple_closurek} can be
simplified as given below
\begin{equation}
[AS_kI]\approx \frac{k-1}{k} \frac{[AS_k][S_kI]}{[S_k]} \approx \frac{k-1}{k} \frac{[AS][SI]k^2[S_k]}{S_1^2} = \frac{[AS][SI]k(k-1)[S_k]}{S_1^2} ,
\label{cPWclosure}
\end{equation}
where $S_1=\sum\limits_{k=1}^N k[S_k]$ is the first moment of the distribution of susceptible nodes. This leads to the so-called compact pairwise model, in which the variables are: $[SI]$, $[SS]$, $[II]$ and $[S_k]$ for $k=1,2,\ldots , K$, i.e., it contains $K+3$ differential equations. In fact, the system consists of equations \eqref{SIhetde1}, and \eqref{SIpairde3}-\eqref{SIpairde5} with the above mentioned closures and approximations, namely
\eqref{cPWappr} and \eqref{cPWclosure}. Thus it takes the form
\begin{eqnarray}
\dot {[S_k]_c} &=& \gamma [I_k]_c - \tau k [S_k]_c \frac{[SI]_c}{S_s}, \label{cPWde1}\\
\dot {[I_k]_c} &=& \tau k [S_k]_c \frac{[SI]_c}{S_s} - \gamma [I_k]_c, \label{cPWde2}\\
\dot {[SI]_c} &=& \gamma ([II]_c-[SI]_c) + \tau ([SS]_c-[SI]_c)[SI]_c P - \tau [SI]_c, \label{cPWde3}\\
\dot {[SS]_c} &=& 2\gamma [SI]_c - 2\tau [SS]_c[SI]_c P  \label{cPWde4},\\
\dot {[II]_c} &=& 2\tau [SI]_c -2\gamma [II]_c + 2\tau [SI]_c^2 P , \label{cPWde5}\\
S_s &=& \sum_{k=1}^K k[S_k]_c , \quad P= \frac{1}{S_s^2} \sum\limits_{k=1}^K (k-1)k[S_k]_c . \label{cPWde6}
\end{eqnarray}
Here the subscript $c$, referring to the word `compact', is used to emphasize that the solution of this system is different from the exact expected values.

The next level of complexity is represented by the pre-compact pairwise model, in which the variables $[S_kI]$ and $[S_kS]$ are kept as independent variables and differential equations for these are written down. Thus the systems can be formulated in terms of variables such as, $[S_k]$, $[S_kS]$, $[S_kI]$, $[I_kS]$ and $[I_kI]$, i.e., resulting in a total of $5K$ variables. This can be done by considering the closure introduced in \cite{eames2002modeling} which is
\begin{equation}
[A_nB_m]=\frac{[A_nB][A_nB]}{[AB]}\frac{[N_nN_m]\sum_{q}q[N_q]}{n[N_n]m[N_m]}, \label{Eames_Pair_Closure}
\end{equation}
where $N_k$ denotes the number of nodes of degree $k$. It is wroth noting that this system is not able to account for preferential mixing.

The most complex system, which we call heterogeneous pairwise model, uses all combinations of pairs as variables, namely $[S_kS_l]$, $[S_kI_l]$ and $[I_kI_l]$. Hence, it consists of $2K^2$ differential equations. At the price of having a system with the number of equations of quadratic order, we do not need any extra approximations (besides the closures), such as \eqref{cPWappr} in the compact pairwise model, or \eqref{Eames_Pair_Closure} for the pre-compact pairwise model. Without explicitly including the closures, the most complex system can be written as
\begin{eqnarray}
\dot {[S_k]} &=& - \tau \sum{l}[S_kI_l]+\gamma [I_k] , \label{SI_het_exact_Sk}\\
\dot {[I_k]} &=& +\tau \sum{l}[S_kI_l]-\gamma [I_k] , \label{SI_het_exact_Ik}\\
\dot {[S_kS_l]} &=& -\tau \sum_{m}\left([I_mS_kS_l]+[S_kS_lI_m]\right)+\gamma \left([S_kI_l]+[I_kS_l]\right), \label{SI_het_exact_SkSl}\\
\dot {[S_kI_l]} &=& +\tau \sum_{m}\left([S_kS_lI_m]-[I_mS_kI_l]\right)-(\tau+\gamma) [S_kI_l]+\gamma [I_kI_l], \label{SI_het_exact_SkIl}\\
\dot {[I_kI_l]} &=& +\tau \sum_{m}\left([I_mS_kI_l]+[I_kS_lI_m]\right)+\tau \left( [S_kI_l]+[I_kS_l] \right)-2\gamma [I_kI_l], \label{SI_het_exact_IkIl}
\end{eqnarray}
with all subscripts going from $1,2, \dots, K$.

\subsection{Super compact pairwise model with heterogeneous triple closure}

We now show that the network heterogeneity can be captured by a small system, containing only four differential equations, just as in the simplest pairwise model. Consider a random network with degrees $d_1, d_2, \ldots , d_K$ and denote the number of nodes of degree $d_k$ by $N_k$ for $k=1,2, \ldots , K$, i.e., $N_1+N_2+\ldots
+ N_K=N$. We note that denoting degrees as $d_k$ instead of $k$ will prove to be advantageous in the derivation below.
The degree distribution of the graph is then given by $p_k=\frac{N_k}{N}$. The average degree and the second moment of the degree distribution are
\begin{equation}
\langle k \rangle = \frac{1}{N}\sum\limits_{k=1}^K d_k N_k  , \quad \langle k^2 \rangle = \frac{1}{N}\sum\limits_{k=1}^K d_k^2 N_k .
\label{avkk2}
\end{equation}

In order to arrive to our new even more simplified system, the super compact PW model, we start from a triple and the closure given in \eqref{cPWclosure}
$$
[ASI] = \sum_{k=1}^N [AS_kI] \approx \frac{[AS][SI]}{S_1^2}  \sum_{k=1}^N d_k(d_k-1)[S_k] = [AS][SI] \frac{S_2-S_1}{S_1^2},
$$
where we used closures \eqref{triple_closurek} and \eqref{cPWappr}, and where $S_2=\sum\limits_{k=1}^K d_k^2[S_k]$ is the second moment of the distribution of susceptible nodes. Thus in order to use this closure in the exact system \eqref{SIpairde1}-\eqref{SIpairde5} one needs an algebraic expression of $S_2$ and $S_1$ in terms of variables $[S]$, $[I]$, $[SI]$, $[II]$ and $[SS]$ only. Expressing the total number of stubs starting from susceptible nodes we get $S_1=[SI] + [SS]$ as an exact relation. Thus the problem arises from the fact that such an exact relation is not available for the second moment $S_2$. Our heuristic idea in obtaining a good approximation of $\frac{S_2-S_1}{S_1^2}$ is the following. Dividing the equation $[S]=\sum\limits_{k=1}^K [S_k]$ by $[S]$ we get that $[S_k]/[S]$ is a probability distribution. The expected value of this distribution is known, it is
$$
\sum\limits_{k=1}^K d_k \frac{[S_k]}{[S]} = n_S := \frac{[SI] + [SS]}{[S]},
$$
or in other words the average degree of susceptible nodes. Our idea is to use an a priori approximating distribution for $[S_k]/[S]$ that will be denoted by $s_k$. This approximating distribution satisfies
\begin{eqnarray}
s_1+s_2+ \ldots + s_K &=& 1, \label{sk1} \\
d_1s_1+d_2s_2+ \ldots + d_Ks_K &=& n_S . \label{sk2}
\end{eqnarray}
In order to get an a priori approximating distribution we determined $[S_k]/[S]$ numerically from the compact pairwise model and compared it to $p_k=N_k/N$, the degree distribution of the graph. Numerical results show that these are linearly related, meaning that $s_k/p_k$ is a linear function of the degree $d_k$. More precisely,  $s_k/p_k$ can be written as $A(t)d_k+B(t)$, where $A$ and $B$ are time dependent with this relation assumed to hold for all degrees. This allows to deal with the heavily under determined linear system given by Eqs. (\ref{sk1})-(\ref{sk2}). Introducing the notation $q_k=s_k/p_k$ the assumption on linearity can be formulated as
$$
\frac{q_k-q_1}{d_k-d_1} = \frac{q_K-q_1}{d_K-d_1} , \quad k=1,2, \ldots , K.
$$
This yields an expression for $q_k$ in terms of $q_1, q_K$ and the degrees $d_k$ as
$$
(d_K-d_1) q_k = (d_k-d_1) q_K + (d_K-d_k) q_1 .
$$
Multiplying this equation by $p_k$ we get the following relation between $s_k$ and $p_k$
\begin{equation}
(d_K-d_1) s_k = p_k(d_k-d_1) q_K + p_k(d_K-d_k) q_1 .
\label{sk}
\end{equation}
Observe that $q_1$ and $q_K$ can be determined from system \eqref{sk1}-\eqref{sk2} by substituting the above expression for $s_k$. Namely, we obtain
\begin{eqnarray}
(d_K-d_1) &=& (n_1-d_1)q_K + (d_K-n_1) q_1, \label{qeq1} \\
(d_K-d_1) n_S &=& (n_2-n_1 d_1)q_K + (n_1 d_K-n_2) q_1, \label{qeq2}
\end{eqnarray}
where $n_i=\sum\limits_{k=1}^K d_k^i p_k$ is the $i$-th moment of the degree distribution. (It is more convenient to use $n_1$ and $n_2$ instead of $\langle k \rangle$ and $\langle k^2 \rangle$.) Solving the linear system \eqref{qeq1}-\eqref{qeq2} for $q_1$ and $q_K$ we get
\begin{eqnarray}
(n_2-n_1^2)q_1 &=& n_2-n_1n_S + d_1(n_S-n_1), \label{q1} \\
(n_2-n_1^2)q_K &=& n_2-n_1n_S + d_K(n_S-n_1). \label{q2}
\end{eqnarray}
Substituting these expressions into \eqref{sk} leads to
$$
(d_K-d_1) (n_2-n_1^2) s_k = p_k(d_k-d_1) (n_2-n_1n_S + d_K(n_S-n_1)) + p_k(d_K-d_k) (n_2-n_1n_S + d_1(n_S-n_1)) .
$$
Now we are in the position of determining the approximate second moment of the distribution $s_k$. Multiplying the above equation by $d_k^2$ and summing from $k=1$ to $k=K$ some simple algebra yields
$$
(n_2-n_1^2) \sum\limits_{k=1}^K d_k^2 s_k = n_2(n_2-n_Sn_1)+ n_3(n_S-n) .
$$
Note that the third moment $n_3$ of the degree distribution comes into play. Thus the desired quantity $S_2$ can be approximated as
$$
S_2=\sum\limits_{k=1}^K d_k^2[S_k] \approx \sum\limits_{k=1}^K d_k^2 [S] s_k = [S] \frac{n_2(n_2-n_Sn_1)+ n_3(n_S-n)}{n_2-n_1^2}.
$$
Hence using $S_1=[SI] + [SS] = n_S [S]$ we get
$$
\frac{S_2-S_1}{S_1^2} \approx  \frac{1}{n_S^2 [S]} \left( \frac{n_2(n_2-n_Sn_1)+ n_3(n_S-n)}{n_2-n_1^2} -n_S \right).
$$
Therefore, the new closure relation is
\begin{equation}
[ASI] = \frac{[AS][SI]}{n_S [S]} \left( \frac{n_2(n_2-n_Sn_1)+ n_3(n_S-n)}{n_S(n_2-n_1^2)} -1 \right) .
\label{newcl}
\end{equation}

We note that in the case of a homogeneous network, where each node has degree $n$, we have $n_S=n$ and the average degree is $n_1=n$. Hence, the expression in the bracket simplifies to $\frac{n_2}{n} -1$. Moreover, the second moment is $n_2 = n^2$. Therefore, this term is simply ($n-1$) and leads to the traditional closure $[ASI] = \frac{n-1}{n} \frac{[AS][SI]}{[S]}$.

Using the new closure \eqref{newcl} in the pairwise model \eqref{SIpairde1}-\eqref{SIpairde5} gives the super compact PW model in the following form:
\begin{eqnarray}
\dot {[S]_s} &=& \gamma [I]_s - \tau [SI]_s, \label{scPW1}\\
\dot {[I]_s} &=& \tau [SI]_s- \gamma [I]_s, \label{scPW2}\\
\dot {[SI]_s} &=& \gamma ([II]_s-[SI]_s) + \tau [SI]_s ([SS]_s-[SI]_s)Q -\tau [SI]_s, \label{scPW3}\\
\dot {[SS]_s} &=& 2\gamma [SI]_s -2 \tau [SI]_s [SS]_s Q, \label{scPW4} \\
\dot {[II]_s} &=& -2\gamma [II]_s + 2\tau [SI]_s^2 Q+ 2\tau [SI]_s , \label{scPW5}
\end{eqnarray}
where
$$
Q=\frac{1}{n_S [S]} \left( \frac{n_2(n_2-n_Sn_1)+ n_3(n_S-n)}{n_S(n_2-n_1^2)} -1 \right), \quad  n_S := \frac{[SI] + [SS]}{[S]}.
$$
In the next section we show that this new super compact pairwise model gives an as accurate output as the compact pairwise model, despite of the fact that it
contains significantly fewer differential equations.

\section{Performance of the new closure for different networks}

As it was shown in the Introduction in Figure \ref{fig:fig1}, the heterogeneous PW and compact PW models give very good agreement with simulations, hence we compare the super compact PW model with the new closure to the compact PW model. This comparison will be done for different heterogeneous networks. Thus systems \eqref{SIpairclosde1}-\eqref{SIpairclosde5}, \eqref{cPWde1}-\eqref{cPWde5} and \eqref{scPW1}-\eqref{scPW5} will be solved numerically and the time dependence of $[I]_p$, $[I]_c$ and $[I]_s$ are compared, where $[I]_c = \sum\limits_{k=1}^N [I_k]_c$ is the total number of infected nodes in the compact PW model. The parameters of the epidemic are fixed at $\gamma =1$ and $\tau = 3 \gamma \langle k \rangle/\langle k^2 \rangle$. The later is chosen in such a way that the ratio of $\tau$ and its critical value $\tau_{cr} = \gamma \langle k \rangle/\langle k^2 \rangle$ is a given constant. Here, this ratio is chosen to be 3, its actual value has only a minor influence on the results, generally this need to be greater than 1 to have an epidemic.

\begin{figure}[h!]
  \centering
	\includegraphics[width=0.6\textwidth, height=0.3\textheight]{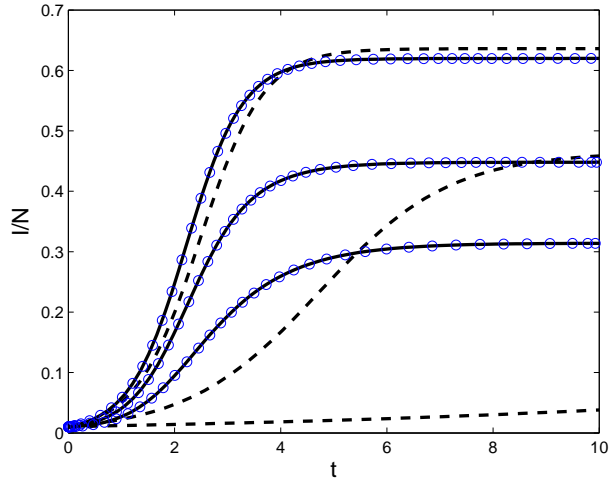}
	\caption{The curves $[I]_p$ (dashed), $[I]_c$ (continuous) and $[I]_s$ (circles) for a bimodal graph with different ratios of the number of low and high degree nodes. The upper curves correspond to $N_1=0.1N$, $N_2=0.9N$, the middle ones are based on $N_1=0.5N$, $N_2=0.5N$ and the lower are for $N_1=0.9N$, $N_2=0.1N$. The parameter values are $N=1000$, $k_1=5$, $k_2=35$, $\gamma=1$ and $\tau = 3 \gamma \langle k \rangle/\langle k^2 \rangle$.} \label{fig:bimod}
\end{figure}

Let us consider first the case of bimodal random graphs, where there are two different degrees $d_1$ and $d_2$, $N_1$ denotes the number of nodes with degree $d_1$ and $N_2$ denote the number of nodes with degree $d_2$, that is $N_1+N_2=N$. In order to investigate the effect of graph structure the ratio of low and high degree nodes, i.e., $N_1$ and $N_2$ is varied. The degrees are fixed at $k_1=5$ and $k_2=35$. In Figure \ref{fig:bimod} the curves $[I]_p$, $[I]_c$ and $[I]_s$ are shown in three cases. The average degree and the standard deviation of the degree distribution is shown in Table \ref{table:q}.
One can see that the new system agrees with and is almost indistinguishable from the compact pairwise model, in fact for bimodal graphs $[I]_s$ coincides with $[I]_c$ since Eqs. (\ref{sk1})-(\ref{sk2}) provide a unique solution without involving any approximations.  Figure \ref{fig:bimod} shows that the traditional pairwise model performs relatively well only in the case when the standard deviation is small, that is the graph is nearly homogeneous.

\begin{figure}[h!]
  \centering
	\includegraphics[width=0.6\textwidth, height=0.3\textheight]{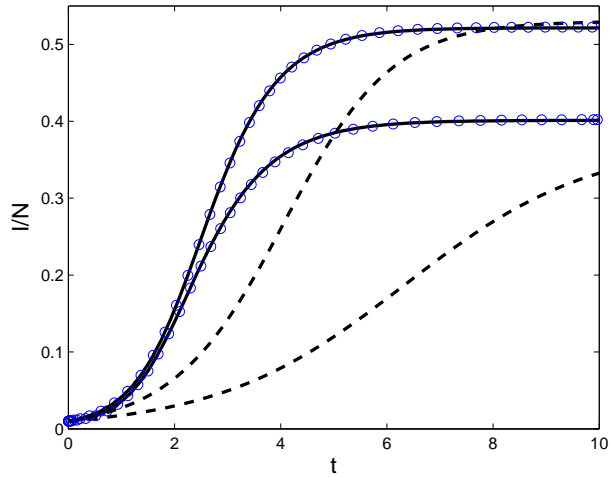}
	\caption{The curves $[I]_p$ (dashed), $[I]_c$ (continuous) and $[I]_s$ (circles) for sparse (lower curves) and a dense (upper curves) power law configuration graphs. The lower curves belong to the sparse case with $k_{min}=5$ and $k_{max}=30$. The upper curves belong to the dense case with  $k_{min}=10$ and $k_{max}=140$. The power is $\alpha = 2$ in both cases. The parameter values are $N=1000$, $\gamma=1$ and $\tau = 3 \gamma \langle k \rangle/\langle k^2 \rangle$.} \label{fig:PL}
\end{figure}

Consider now the case of configuration random graphs with cutoff power law degree distribution. These random graphs are given by a minimal degree $k_{min}$, a maximal degree $k_{max}$ and a power $\alpha$. The degree distribution of the graph is $p(k)=Ck^{-\alpha}$ for $k=k_{min}, k_{min}+1, \ldots , k_{max}$ with the normalisation constant $C$ given by $$\frac{1}{C}=\sum\limits_{k=k_{min}}^{k_{max}} k^{-\alpha}.$$
In Figure \ref{fig:PL} again the functions $[I]_p$, $[I]_c(t)$ and $[I]_s(t)$ are shown for a sparse (lower curves) and a dense (upper curves) power law configuration graph with power $\alpha = 2$. Table \ref{table:q} again shows the average degree and the standard deviation of the degree distribution of the sparse and dense networks.
The value of $\tau$ in both cases is $\tau = 3 \gamma \langle k \rangle/\langle k^2 \rangle$. We can see again that the super compact PW model gives excellent agreement with the compact pairwise model.

\begin{center}
\begin{table}[ht]
\begin{tabular}{ | c | c | c | }
\hline
Network  & $\langle k \rangle$ & $\sqrt{\langle k^2 \rangle-\langle k \rangle^2}$  \\
  \hline
  Bimodal 0.1 & 32 & 9  \\
  \hline
  Bimodal 0.5 & 20 & 15  \\
  \hline
  Bimodal 0.9 & 8 & 9  \\
  \hline
  Power law sparse & 10.1 & 5.9  \\
  \hline
  Power law dense & 28.4 & 26.01 \\
  \hline
\end{tabular}
\caption{The average degree and the standard deviation of the degree distribution of the graphs for which the performance of the new closure was tested. For bimodal graphs the degrees are $k_1=5$ and $k_2=35$, the numbers in the first coloumn indicate the proportion of low degree nodes, i.e., $N_1/N$. For the sparse power law graphs the degrees vary between $k_{min}=1$ and $k_{max}=35$, for the dense one $k_{min}=10$ and $k_{max}=140$, the power is $\alpha = 2$. }
\label{table:q}
\end{table}
\end{center}

\section{Discussion}

In this short paper, we derived a super compact pairwise model consisting of only 4 equations for $SIS$ dynamics and for heterogenous networks constructed according to the
configuration model. This represents an improvement of going from order $K$, where $K$ is the number of distinct degrees in the network, to order one equations, namely 4. We note
that the closure that made the reduction possible relies on the observation that the distribution of susceptible nodes of degree $k$, which is time dependent, can be related to the original degree distribution of the network via a simple linear relation. We note that the linear relation may not be the single or unique choice, more sophisticated functional forms could be used based on combinatorial arguments. Moreover, the closure will not only encompass the first and second moment of the degree distribution but also the third. The new super compact model gives excellent agreement with the previously derived compact
pairwise model.

The accuracy of the new closure can be estimated in a semi-analytic way. The numerical solution of the compact PW will allow to evaluate
$$
E=\frac{S_2-S_1}{S_1^2} - Q,
$$
which quantifies the performance of the newly derived closure, upon using the compact PW model as a benchmark. Moreover, it can be shown analytically
that the difference $|[I]_s(t) - [I]_c(t)|$ can be estimated by a constant multiple of $E$. Further work on this model will include a detailed bifurcation study of the closed super compact pairwise system and we will aim to determine the stability of the disease free and endemic steady states. We will also investigate whether an explicit formula for the endemic steady state is possible. If these calculations will be tractable, the stability of the disease free steady state should also yield $R_0$ or at least an $R_0$-like quantity for the super compact PW model.
%
%

\bibliography{Sup_Comp_Pairwise}
\bibliographystyle{ieeetr}

\end{document}